\documentclass[11pt]{article}
\usepackage[margin=1in]{geometry} 
\usepackage{amssymb,amsfonts,amsmath,bbm,mathrsfs,stmaryrd,mathtools}
\usepackage{xcolor}
\usepackage{url}

\usepackage[shortlabels]{enumitem}
\usepackage{tensor}

\usepackage{xr}
\usepackage[colorlinks,
linkcolor=black!75!red,
citecolor=blue,
pdftitle={},
pdfproducer={pdfLaTeX},
pdfpagemode=None,
bookmarksopen=true,
bookmarksnumbered=true,
backref=page]{hyperref}

\usepackage{tikz}
\usetikzlibrary{arrows,calc,decorations.pathreplacing,decorations.markings,intersections,shapes.geometric,through,fit,shapes.symbols,positioning,decorations.pathmorphing}

\usepackage{braket}

\usepackage[amsmath,thmmarks,hyperref]{ntheorem}
\usepackage{cleveref}

\creflabelformat{enumi}{#2#1#3}

\crefname{section}{Section}{Sections}
\crefformat{section}{#2Section~#1#3} 
\Crefformat{section}{#2Section~#1#3} 

\crefname{subsection}{\S}{\S\S}
\AtBeginDocument{%
  \crefformat{subsection}{#2\S#1#3}%
  \Crefformat{subsection}{#2\S#1#3}%
}

\crefname{subsubsection}{\S}{\S\S}
\AtBeginDocument{%
  \crefformat{subsubsection}{#2\S#1#3}%
  \Crefformat{subsubsection}{#2\S#1#3}%
}

\theoremstyle{plain}

\newtheorem{lemma}{Lemma}[section]

\newtheorem{corollary}[lemma]{Corollary}
\newtheorem{theorem}[lemma]{Theorem}

\theoremstyle{plain}
\theoremnumbering{Alph}

\theoremstyle{plain}
\theorembodyfont{\upshape}
\theoremsymbol{\ensuremath{\blacklozenge}}

\newtheorem{remark}[lemma]{Remark}

\crefname{definition}{definition}{definitions}
\crefformat{definition}{#2definition~#1#3} 
\Crefformat{definition}{#2Definition~#1#3} 

\crefname{ex}{example}{examples}
\crefformat{example}{#2example~#1#3} 
\Crefformat{example}{#2Example~#1#3} 

\crefname{exs}{example}{examples}
\crefformat{examples}{#2example~#1#3} 
\Crefformat{examples}{#2Example~#1#3} 

\crefname{remark}{remark}{remarks}
\crefformat{remark}{#2remark~#1#3} 
\Crefformat{remark}{#2Remark~#1#3} 

\crefname{remarks}{remark}{remarks}
\crefformat{remarks}{#2remark~#1#3} 
\Crefformat{remarks}{#2Remark~#1#3} 

\crefname{convention}{convention}{conventions}
\crefformat{convention}{#2convention~#1#3} 
\Crefformat{convention}{#2Convention~#1#3} 

\crefname{notation}{notation}{notations}
\crefformat{notation}{#2notation~#1#3} 
\Crefformat{notation}{#2Notation~#1#3} 

\crefname{table}{table}{tables}
\crefformat{table}{#2table~#1#3} 
\Crefformat{table}{#2Table~#1#3}

\crefname{lemma}{lemma}{lemmas}
\crefformat{lemma}{#2lemma~#1#3} 
\Crefformat{lemma}{#2Lemma~#1#3} 

\crefname{proposition}{proposition}{propositions}
\crefformat{proposition}{#2proposition~#1#3} 
\Crefformat{proposition}{#2Proposition~#1#3} 

\crefname{propositionN}{proposition}{propositions}
\crefformat{propositionN}{#2proposition~#1#3} 
\Crefformat{propositionN}{#2Proposition~#1#3} 

\crefname{corollary}{corollary}{corollaries}
\crefformat{corollary}{#2corollary~#1#3} 
\Crefformat{corollary}{#2Corollary~#1#3} 

\crefname{corollaryN}{corollary}{corollaries}
\crefformat{corollaryN}{#2corollary~#1#3} 
\Crefformat{corollaryN}{#2Corollary~#1#3} 

\crefname{theorem}{theorem}{theorems}
\crefformat{theorem}{#2theorem~#1#3} 
\Crefformat{theorem}{#2Theorem~#1#3} 

\crefname{theoremN}{theorem}{theorems}
\crefformat{theoremN}{#2theorem~#1#3} 
\Crefformat{theoremN}{#2Theorem~#1#3} 

\crefname{enumi}{}{}
\crefformat{enumi}{#2#1#3}
\Crefformat{enumi}{#2#1#3}

\crefname{assumption}{assumption}{Assumptions}
\crefformat{assumption}{#2assumption~#1#3} 
\Crefformat{assumption}{#2Assumption~#1#3} 

\crefname{construction}{construction}{Constructions}
\crefformat{construction}{#2construction~#1#3} 
\Crefformat{construction}{#2Construction~#1#3} 

\crefname{equation}{}{}
\crefformat{equation}{(#2#1#3)} 
\Crefformat{equation}{(#2#1#3)}

\numberwithin{equation}{section}

\theoremstyle{nonumberplain}
\theoremsymbol{\ensuremath{\blacksquare}}

\newcommand\pf[1]{\newtheorem{#1}{Proof of \Cref{#1}}}

\newcommand\bG{{\mathbb G}}
\newcommand\bH{{\mathbb H}}

\newcommand\bU{{\mathbb U}}

\newcommand\bZ{{\mathbb Z}}

\newcommand\fu{{\mathfrak u}}

\DeclareMathOperator{\Ad}{Ad}

\newcommand{\cat}[1]{\textsc{#1}}

\newcommand{\qedhere}{\mbox{}\hfill\ensuremath{\blacksquare}}

\title{Straightening Banach-Lie-group-valued almost-cocycles}
\author{Alexandru Chirvasitu}

\begin{document}

\date{}

\newcommand{\Addresses}{{
  \bigskip
  \footnotesize

  \textsc{Department of Mathematics, University at Buffalo}
  \par\nopagebreak
  \textsc{Buffalo, NY 14260-2900, USA}  
  \par\nopagebreak
  \textit{E-mail address}: \texttt{achirvas@buffalo.edu}


}}

\maketitle

\begin{abstract}
  For a compact group $\mathbb{G}$ acting continuously on a Banach Lie group $\mathbb{U}$, we prove that maps $\mathbb{G}\to \mathbb{U}$ close to being 1-cocycles for the action can be deformed analytically into actual 1-cocycles. This recovers Hyers-Ulam stability results of Grove-Karcher-Ruh (trivial $\mathbb{G}$-action, compact Lie $\mathbb{G}$ and $\mathbb{U}$) and de la Harpe-Karoubi (trivial $\mathbb{G}$-action, $\mathbb{U}:=$invertible elements of a Banach algebra). The obvious analogues for higher cocycles also hold for abelian $\mathbb{U}$.
\end{abstract}


\vspace{.5cm}



\section*{Introduction}

The motivation for the present note was the juxtaposition in \cite[\S 5]{wein_lin} of the following two results, referred to there as the {\it almost homomorphism} and {\it almost representation} theorems respectively:
\begin{itemize}
\item Continuous maps $\bG\xrightarrow{\rho} \bH$ between compact Lie groups that are {\it almost} morphisms in the sense that
  \begin{equation*}
    \sup_{s,t\in \bG} d(\rho(s)\rho(t),\ \rho(st))\text{ is small}
  \end{equation*}
  for an appropriate metric $d$ on $\bH$ are uniformly close to actual morphisms \cite[Theorem 4.3]{gkr_actcurv_invent_1974}. 
\item Similarly, continuous almost-morphisms
  \begin{equation*}
    \bG
    \xrightarrow{\quad}
    A^{\times}:=\text{invertible elements in a Banach algebra }A
  \end{equation*}
  are close to morphisms \cite[Proposition 4]{zbMATH03577502}. 
\end{itemize}

These are instances of {\it Hyers-Ulam stability} (\cite{zbMATH00149467} and its references, say), and admit a common generalization (stated formally in \Cref{cor:almost-cpct2ban} below): continuous almost-morphisms from a compact group $\bG$ to a {\it Banach Lie group} \cite[Definition IV.I]{neeb-inf} $\bU$ are uniformly close to morphisms. The two results can then be recovered by specializing the Banach Lie side of the picture:  
\begin{itemize}
\item $\bU$ compact (and hence \cite[Theorem IV.3.16]{neeb-lc} finite-dimensional) yields \cite[Theorem 4.3]{gkr_actcurv_invent_1974};
  
\item while $\bU:=A^{\times}$ for a Banach algebra $A$ returns \cite[Proposition 4]{zbMATH03577502}. 
\end{itemize}

The argument has a {\it Newton-approximation} \cite[Definition 1.6]{ms_numan} flavor, very much in the spirit of the one delivering \cite[Theorem 3.1]{john_approx} (say). That argument suggests (and very little additional effort yields) natural generalizations: morphisms $\bG\to \bU$ are nothing but {\it $\bU$-valued $1$-cocycles} on $\bG$ \cite[Appendix to Chapter VII, p.123]{ser_locf} for the trivial $\bG$-action on $\bU$. Dropping that triviality (\Cref{th:almost-cpct2ban}):
\begin{itemize}
\item Almost-1-cocycles $\bG\to \bU$ with respect to a fixed action of a compact group $\bG$ on a Banach Lie group $\bU$ are close to actual 1-cocycles.

\item The same goes for $n$-cocycles with arbitrary $n\in \bZ_{>0}$ provided $\bU$ is abelian. 
\end{itemize}

\subsection*{Acknowledgements}

I am grateful for helpful comments from K.-H. Neeb.

This work is partially supported by NSF grant DMS-2001128. 

\section{Preliminaries}\label{se:prel}

We follow \cite[\S\S 2.3, 3.2, 5.1 and 5.3]{bourb_vars_1-7} (and \cite[Chapter III, introductory remarks]{bourb_lie_1-3}) in adopting a uniform notation for $C^r$ (Banach) manifolds and morphisms for
\begin{equation*}
  r\in \widetilde{\bZ}_{\ge 0}:=\overline{\bZ}_{\ge 0}\sqcup\{\omega\}
  ,\quad
  \overline{\bZ}_{\ge 0}:=\bZ_{\ge 0}\sqcup\{\infty\}
  ,\quad
  \bZ_{\ge 0}\ni n<\infty<\omega,
\end{equation*}
where $C^r$, $r\in \overline{\bZ}_{\ge 0}$ means, as usual, {\it $r$-times continuously differentiable} and $C^{\omega}$ means {\it (real or complex) analytic}.

Write $C^r(X,Y)$ for $C^r$ maps between $C^{\ge r}$ manifolds $X$ and $Y$, with $X$ typically finite-dimensional if $r\ge 1$; we will often abbreviate $C^0(-,-)$ to $C(-,-)$. The {\it $C^r$ topology} \cite[\S 10, Example I]{trev_tvs} on $C^r(-,-)$ is that of uniform convergence on compact sets for all derivatives up to order $r$. In particular, the $C^{\infty}$ and $C^{\omega}$ topologies coincide. 

Consider a compact group $\bG$ and a Banach Lie group $\bU$, the former acting continuously on the latter by automorphisms via
\begin{equation*}
  \bG\times \bU
  \ni (s,x)
  \xmapsto{\quad}
  \tensor[^s]{x}{}
  \in \bU.
\end{equation*}

We employ the language of group cohomology fairly liberally: cocycles, coboundaries, etc., as can be recalled briefly from \cite[Chapter 9]{rot}, \cite[Chapter VII and its Appendix]{ser_locf}, \cite[Chapter III]{brown_coh-gp}, \cite[Chapter IV]{cf_ant_1967} and many other sources. In particular, given the above setup and $\rho\in C(\bG^n,\ \bU)$, write (e.g. \cite[\S VII.3, equation (**)]{ser_locf})
\begin{equation}\label{eq:cobound}
  \begin{aligned}
    \delta\rho(s_0,\ s_1,\ \cdots,\ s_n)
    :=
    \tensor[^{s_0}]{\rho(s_1,\ \cdots,\  s_n)}{}
    &-\rho(s_0 s_1,\ s_2,\ \cdots,\ s_n)\\
    &+\rho(s_0,s_1 s_2,\ \cdots,\ s_n)\\
    &- \cdots\\
    &+(-1)^{n+1} \rho(s_0,s_1,\ \cdots,\ s_{n-1})
  \end{aligned}
\end{equation}
for the {\it coboundary} of $\rho$. This includes uniformly the cases when $\bU$ is abelian and $n$ is arbitrary, or $\bU$ is arbitrary and $n\le 1$ (per the discussion on non-abelian cohomology in \cite[Appendix to Chapter VII]{ser_locf}), where \Cref{eq:cobound} is interpreted as
\begin{equation*}
  \delta x(s) := \tensor[^s]{x}{}\cdot x^{-1}
  \text{ for }n=0
\end{equation*}
(with $x\in \bU$ interpreted as a function $\bG^0\to \bU$) and
\begin{equation*}
  \delta\rho(s_0,s_1) = \rho(s_0)\cdot \tensor[^{s_0}]{\rho(s_1)}{}\cdot \rho(s_0 s_1)^{-1}
  \text{ for }n=1.
\end{equation*}
Write
\begin{equation*}
  C_Z(\bG^n,\ \bU)
  :=
  \left\{\text{continuous }\bG\xrightarrow{\rho}\bU\ |\ \delta \rho\equiv 1\right\}
\end{equation*}
and similarly with $C^r$ in place of $C:=C^0$ when $\bG$ is Lie. 


\section{Analytic deformations of almost-morphisms}\label{se:definto}

For subsets $S\subseteq \bU$ and $B\subseteq \fu:=Lie(\bU)$ of a Banach Lie group and its Lie algebra respectively we write
\begin{equation*}
  S_{\Ad\subset B}
  :=
  \left\{s\in S\ |\ \Ad(s)\in B\right\}
  ,\quad
  \bU\xrightarrow[]{\Ad:=\text{the {\it adjoint representation} \cite[\S III.3.12]{bourb_lie_1-3}}}GL(\fu)
\end{equation*}
The same notation applies to spaces of maps $X\to \bU$: an $\Ad\subset B$ subscript indicates maps taking values in $B$ after composing with $\Ad$. 

\begin{theorem}\label{th:almost-cpct2ban}
  Let $\bG$ be a compact group, $n\in \bZ_{>0}$, and $\bU$ a Banach Lie group, abelian if $n\ge 2$. Let also $B\subset GL(\fu)$ be a bounded open subset. 
  \begin{enumerate}[(1),wide]
  \item\label{item:th:almost-cpct2ban:c0} For every neighborhood $V\ni 1\in C(\bG^n,\ \bU)$ there is a neighborhood $W\ni 1\in C(\bG^{n+1},\bU)$ and an analytic, uniformly continuous map
    \begin{equation*}
      \left(\delta^{-1}W\right)_{\Ad\subset B}
      \ni \rho
      \xmapsto{\quad}
      \rho'\in C_Z(\bG^n,\ \bU)
      ,\quad
      \rho'\cdot \rho^{-1}\in V\text{ throughout}. 
    \end{equation*}

  \item\label{item:th:almost-cpct2ban:cr} Moreover, if $\bG$ is Lie (and hence also an analytic manifold) the analogous statements hold upon substituting $C^r$ ($1\le r\le \omega$) for $C^0$ throughout. 
  \end{enumerate}
\end{theorem}

We first record the consequence mentioned in the introduction. Specializing to $n=1$ and a a {\it trivial} $\bG$-action on $\bU$, 1-cocycles are nothing but topological-group morphisms $\bG\to \bU$; thus:

\begin{corollary}\label{cor:almost-cpct2ban}
  Let $\bG$ be a compact group, $\bU$ a Banach Lie group, and $B\subset GL(\fu)$ be a bounded open subset. 
  \begin{enumerate}[(1),wide]
  \item For every neighborhood $V\ni 1\in C(\bG,\ \bU)$ there is a neighborhood $W\ni 1\in C(\bG^{2},\bU)$ and an analytic, uniformly continuous map
    \begin{equation*}
      \left(\delta^{-1}W\right)_{\Ad\subset B}
      \ni \rho
      \xmapsto{\quad}
      \rho'\in \cat{TopGp}(\bG,\bU)
      ,\quad
      \rho'\cdot \rho^{-1}\in V\text{ throughout}. 
    \end{equation*}

  \item If $\bG$ is Lie the analogous statements hold upon substituting $C^r$ ($1\le r\le \omega$) for $C^0$.  \qedhere
  \end{enumerate}
\end{corollary}

Closeness estimates in the Lie algebra $\fu:=Lie(\bU)$ are with respect to a complete norm $\|\cdot\|$ thereon with
\begin{equation*}
  \|[x,y]\| \le C\|x\|\cdot \|y\|
  ,\quad\forall x,y\in \fu
  \quad\text{for some }C>0;
\end{equation*}
such a norm always exists \cite[\S III.3.7, Corollary to Proposition 24]{bourb_lie_1-3}, and we can of course always scale to $C=1$. $\bU$ is itself completely metrizable \cite[\S III.1.1, Proposition 1]{bourb_lie_1-3}, and the {\it exponential map} \cite[\S III.4.3, Theorem 4]{bourb_lie_1-3}
\begin{equation*}
  \fu
  \xrightarrow{\quad\exp=e^{(\cdot)}\quad}
  \bU
\end{equation*}
implements an analytic isomorphism between origin neighborhoods of $\fu$ and $\bU$. Closeness between ``small'' elements thus transports over back and forth between the two spaces. 

\pf{th:almost-cpct2ban}
\begin{th:almost-cpct2ban}
  We focus on the case $n=1$; the abelianness of $\bU$ makes the argument, if anything, even simpler otherwise. Moreover, the analytic nature of the construction $\rho\xmapsto \rho'$ will make it clear that that construction preserves $r$-fold continuous differentiability or indeed analyticity, so that part \Cref{item:th:almost-cpct2ban:cr} need not be addressed separately.
  
  The assumption is that
  \begin{equation}\label{eq:rhostbeta}
    \rho(s)\cdot \tensor[^s]{\rho(t)}{} = e^{\beta(s,t)}\rho(st),\quad s,t\in \bG
  \end{equation}
  with
  \begin{equation*}
    \left(
      \bG^2\ni (s,t)
      \xmapsto{\quad}
      \beta(s,t)
      \in
      \fu
    \right)
    =
    O(\varepsilon)
    \text{ for small }\varepsilon>0
  \end{equation*}
  in standard {\it big-oh notation} \cite[\S 3.2]{clrs_alg-4e}:
  \begin{equation}\label{eq:bigohk}
    \|\beta(-,-)\|\le K\varepsilon
    \text{ uniformly for some }K>0,
  \end{equation}
  valid universally, for any $\varepsilon>0$, provided the latter is sufficiently small. 
  
  The goal is to produce 
  \begin{equation*}
    \left(
      \bG\ni s
      \xmapsto{\quad}
      \alpha(s)\in \fu
    \right)
    =O(\varepsilon)
  \end{equation*}
  again, attached analytically to $\rho$, so that 
  \begin{equation}\label{eq:targeteq}
    \begin{aligned}
      e^{\alpha(s)}\cdot e^{s\triangleright \alpha(t)}
      \cdot \rho(s)\cdot \tensor[^s]{\rho(t)}{}
      &= e^{\alpha(s)}\rho(s)\cdot
        \tensor[^s]{\left(e^{\alpha(t)}\rho(t)\right)}{}\\
      &= e^{\alpha(st)}\rho(st)\\
      &= e^{\alpha(st)} \cdot e^{-\beta(s,t)} \cdot \rho(s)\rho(t)
        ,\quad s,t\in \bG
    \end{aligned}           
  \end{equation}
  meaning that
  \begin{equation*}
    e^{\alpha(s)}\cdot e^{s\triangleright \alpha(t)}
    =
    e^{\alpha(st)} \cdot e^{-\beta(s,t)}
    ,\quad s,t\in \bG,
  \end{equation*}
  where the $\varepsilon$-almost-action `$\triangleright$' is defined by
  \begin{equation*}
    e^{s\triangleright x}
    :=
    \rho(s)\cdot e^x\cdot \rho(s)^{-1}
    ,\quad
    s\in \bG,\ x\in \fu.
  \end{equation*}

  As hinted above, we will construct $\alpha$ by successive approximation as
  \begin{equation}\label{eq:alphaseries}
    \alpha = \log\left(\cdots e^{\alpha_2}\cdot e^{\alpha_1}\right)
    ,\quad
    \alpha_n = O(\varepsilon^n),
  \end{equation}
  where
  \begin{itemize}
  \item there is a single constant $K>0$ as in \Cref{eq:bigohk}, pertinent to all instances of $O(\varepsilon)$, $O(\varepsilon^2)$, etc., valid throughout the proof;

  \item and the convergence of \Cref{eq:alphaseries} follows from this via the {\it Baker-Campbell-Hausdorff (or BCH) formula} (\cite[Definition IV.1.3]{neeb-lc}, \cite[Part I, \S\S IV.7 and IV.8 and part II, \S V.4]{serre-lag}). 
  \end{itemize}  
  We construct the requisite $\alpha_1=O(\varepsilon)$ first, and then proceed recursively. Repeated application of \Cref{eq:rhostbeta} yields
  \begin{equation*}
    \begin{aligned}
      e^{s\triangleright\beta(t,u)} e^{\beta(s,tu)}\rho(stu)
      &=\rho(s)\cdot e^{\tensor[^s]{\beta(t,u)}{}} \cdot \tensor[^s]{\rho(tu)}{}\\
      &= \rho(s)\cdot
        \tensor[^s]{\left(\rho(t) \cdot \tensor[^{t}]{\rho(u)}{}\right)}{}\\
      &= \rho(s)\cdot \tensor[^s]{\rho(t)}{}\cdot \tensor[^{st}]{\rho(u)}{}\\
      &= e^{\beta(s,t)}\cdot \rho(st)\cdot \tensor[^{st}]{\rho(u)}{}\\
      &= e^{\beta(s,t)}\cdot e^{\beta(st,u)}\cdot \rho(stu),
    \end{aligned}    
  \end{equation*}
  meaning that
  \begin{equation*}
    e^{s\triangleright\beta(t,u)}\cdot e^{\beta(s,tu)}
    =
    e^{\beta(s,t)}\cdot e^{\beta(st,u)}
  \end{equation*}
  and hence $\beta$ is an {\it ($\varepsilon^2$)-almost-cocycle} with respect to `$\triangleright$':
  \begin{equation*}
    \delta_{\triangleright}\beta(s,t,u)
    :=
    s\triangleright\beta(t,u) - \beta(st,u) + \beta(s,tu) - \beta(s,t) = O(\varepsilon^2). 
  \end{equation*}
  The usual (e.g. \cite[I, proof of Theorem 2.3]{moore-ext}) Haar-averaging procedure then also ensures that it is a coboundary of an $O(\varepsilon)$ cochain $\alpha_1$ to order $\varepsilon^2$: setting
  \begin{equation}\label{eq:averaging}
    \alpha_1(t)
    :=
    -\int_{\bG} s\triangleright \beta(s^{-1},t)\ \mathrm{d}\mu_{\bG}(s)
    ,\quad
    \mu_{\bG}:=\text{Haar probability measure on $\bG$},
  \end{equation}
  we have both $\alpha_1=O(\varepsilon)$ and
  \begin{equation*}
    \delta_{\triangleright}\alpha_1 + \beta = O(\varepsilon^2)
    ,\quad
    \delta_{\triangleright}\alpha_1(s,t):=s\triangleright\alpha_1(t)-\alpha_1(st)+\alpha_1(s). 
  \end{equation*}
  This means that with $\alpha_1$ in place of $\alpha$, \Cref{eq:targeteq} holds to order $\varepsilon^2$. We can now proceed recursively, making the substitutions
  \begin{equation*}
    \rho \rightsquigarrow e^{\alpha}\rho
    ,\quad
    \beta \rightsquigarrow \beta + \delta\alpha_1
    ,\quad
    \varepsilon\rightsquigarrow \varepsilon^2.
  \end{equation*}
  The uniformity of the constants featuring implicitly in the $O(\varepsilon^n)$ for varying $n$ follows from the boundedness of $\left\{s\triangleright\right\}_{s\in \bG}$ ensured by restricting attention to $\left(\delta^{-1}W\right)_{\Ad\subset B}$, and the $\alpha_n=O(\varepsilon^n)$ are all constructed analytically (via \Cref{eq:averaging}) in terms of the initial data. 
\end{th:almost-cpct2ban}  
  
\begin{remark}\label{re:jones-almostcocycle}
  We saw in the course of the proof of \Cref{th:almost-cpct2ban} that almost-cocycles are close to coboundaries, so in particular close to actual cocycles. This is presumably the type of result alluded to in passing in \cite[Remark immediately preceding \S 3.4]{jones_fin-act-hyp} (which in turn also refers back to \cite{zbMATH03577502}).
\end{remark}



\addcontentsline{toc}{section}{References}

\Addresses

\end{document}